\documentclass[12pt]{article}

\usepackage{amsthm,color,amsmath,amsfonts,amssymb}

\parskip=\medskipamount
\def\onehalf{{\textstyle\frac12}}
\def\onethird{{\textstyle\frac13}}
\def\threehalf{{\textstyle\frac32}}
\def\onesix{{\textstyle\frac16}}

\def\jet{J^1\tau}
\def\dualJ{J^1\tau^*}

\def\vectorfields#1{{\mathcal X}(#1)}

\def\oneforms#1{{\mathcal X}^*(#1)}
\def\forms#1#2{{\textstyle\bigwedge}^{#1}(#2)}

\def\lie#1{{\mathcal L}_{#1}}
\def\fpd#1#2{\frac{\partial #1}{\partial #2}}
\def\fpdd#1#2#3{\frac{\partial^2 #1}{\partial #2\,\partial{#3}}}
\def\R{\mathbb{R}}
\def\D#1{{\mathcal D}_{#1}}
\def\Dperp#1{{\mathcal D}_{#1}^\perp}

\def\proof{{\sc Proof.}}

\newtheorem{lemma}{\bf Lemma}
\newtheorem{thm}{\bf Theorem}

\def\h#1{{#1}^{h}}

\def\wt#1{\widetilde{#1}}
\def\cR{\wt{R}_{T^{\!*}}}

\def\sp{\mathop{\rm sp}}

\def\hook2{\mathop{\kern -2pt\hbox to 6pt{\hrulefill}
                      \hbox{\vrule\phantom{\vbox to 8pt{}}}}\kern 1pt}

\begin{document}

\title{Lifted tensors and Hamilton-Jacobi separability}

\author{G.\ Waeyaert$^{a}$ and W.\ Sarlet$^{a,b,}\footnote{Corresponding author: willy.sarlet@ugent.be}$ \\
{\small ${}^a$Department of Mathematics, Ghent University }\\
{\small Krijgslaan 281, B-9000 Ghent, Belgium}\\[1mm]
{\small ${}^b$Department of Mathematics and Statistics, La Trobe University}\\
{\small Bundoora, Victoria 3086, Australia}
}
\date{}
\maketitle

\begin{quote}
{\bf Abstract.} {\small Starting from a bundle $\tau:E\rightarrow\R$, the bundle $\pi:\dualJ\rightarrow E$, which is the dual of the first jet bundle $\jet$ and a sub-bundle of $T^*E$, is the appropriate manifold for the geometric description of time-dependent Hamiltonian systems. Based on previous work, we recall properties of the complete lifts of a type $(1,1)$ tensor $R$ on $E$ to both $T^*E$ and $\dualJ$. We discuss how an interplay between both lifted tensors leads to the identification of related distributions on both manifolds. The integrability of these distributions, a coordinate free condition, is shown to produce exactly Forbat's conditions for separability of the time-dependent Hamilton-Jacobi equation in appropriate coordinates. }
\end{quote}

\section{Introduction}

Separation of variables is a classical approach in trying to solve the Hamilton-Jacobi equation for a given Hamiltonian system. There is a vast literature about the subject, dating back to the beginning of the previous century. It is an impossible task to try to give a comprehensive account of this literature, so we will merely list a few references, which have some relevance for the present paper and illustrate that there is still a lively interest in the subject. Benenti has contributed a lot to the evolution of the subject in the past decades, see for example \cite{Ben3}, \cite{Ben} and \cite{Ben2}. A standard work about the classification of coordinate systems in which the Hamilton-Jacobi equation separates for Riemannian spaces of constant curvature is the monograph of Kalnins \cite{Kalnins}. Some other fairly recent contributions with a differential geometric content are \cite{IMM}, \cite{Cra03}, \cite{Cra05}, \cite{BolMat} and \cite{FP}. What all these references have in common is that they discuss aspects of separability of the Hamilton-Jacobi equation for autonomous Hamiltonian systems. In that respect, it is well known that necessary and sufficient conditions for Hamilton-Jacobi separability in the autonomous case were developed by Levi-Civita \cite{L-C}. Much less known is that these conditions were generalized to the case of time-dependent systems by Forbat \cite{Forbat}. Forbat's conditions read
\begin{align}
\fpd{H}{p_i}\left(\fpdd{H}{q^i}{q^j}\fpd{H}{p_j}-\fpdd{H}{q^i}{p_j}\fpd{H}{q^j}\right)& =
\fpd{H}{q^i}\left(\fpdd{H}{p_i}{q^j}\fpd{H}{p_j}-\fpdd{H}{p_i}{p_j}\fpd{H}{q^j}\right),
\label{Forbat1} \\[1mm]
\fpd{H}{p_i}\fpdd{H}{q^i}{t}& = \fpd{H}{q^i}\fpdd{H}{p_i}{t},
\label{Forbat2}
\end{align}
where there is no summation over repeated indices. A study of Hamiltonians of mechanical type which satisfy Forbat's conditions was conducted in \cite{Frans}. A weak point about such conditions is that they can merely test whether the Hamilton-Jacobi equation is separable in the given coordinates: one has to be lucky to have chosen separation coordinates already for the test to give positive results. Our purpose is to develop an intrinsic formulation of Forbat's conditions, i.e.\ to obtain a test for the existence of separation coordinates which in principle can be carried out in any given coordinate chart and should then provide information about the way separation coordinates can be constructed.

There is a variety of possible differential geometric models for time-dependent Hamiltonian systems, see for example \cite{EMR}, \cite{CaRan}, \cite{BEMMR} to cite just a few. As we argued already in \cite{SW}, however, an important point for the model we choose is that it should not incorporate time-dependence via a product structure, because time-dependent coordinate transformations do not preserve such structure. More convenient approaches therefore are those which start from a bundle over $\R$, as is the case, for example, in \cite{MV} and \cite{Sarda}. So briefly, let $E$ be a bundle $\tau:E\rightarrow \R$ with $\dim E=n+1$ and coordinates denoted by $(t,q^i)$. The cotangent bundle $T^*E$, with natural coordinates $(t,q^i,p_0,p_i)$, is often said to be the extended dual $(\jet)^\dagger$ of the first jet bundle $\jet$ of $E$. The quotient bundle $T^*E/\langle dt\rangle$ which we denote by $\dualJ$ then is said to be the dual of $\jet$ (sometimes also called the vertical cotangent bundle). There are natural projections, say $\rho:T^*E\rightarrow \dualJ$ and $\pi:\dualJ\rightarrow E$. Each point $m\in\dualJ$ is an equivalence class of covectors $\langle\alpha\rangle\!\!\!\mod dt$ at $\pi(m)$ and has a well-defined action on vertical tangent vectors to $E$; saying that $m$ has coordinates $(t,q^i,p_i)$ means that $\alpha_{(t,q)}=p_idq^i\!\!\!\mod dt$. A Hamiltonian is a section $h$ of the line bundle $\rho:T^*E\rightarrow\dualJ$. Locally, $h$ defines a function $H$ on $\dualJ$, determined by $h:(t,q,p)\mapsto (t,q,p_0=-H(t,q,p),p)$. If $\omega_E=d\theta_E$ denotes the canonical symplectic form on $T^*E$, we have that locally $h^*\omega_E = dp_i\wedge dq^i - dH\wedge dt$ and $X_h$, defined by $i_{X_h}h^*\omega_E=0$ and $\langle X_h,dt\rangle=1$, is the associated Hamiltonian vector field on $\dualJ$, locally of the form (here with the usual summation convention)
\begin{equation}
X_h = \fpd{}{t} + \fpd{H}{p_i}\fpd{}{q^i} - \fpd{H}{q^i}\fpd{}{p_i}. \label{Xh}
\end{equation}

In principle, one may hope to develop an intrinsic model for the separability issue directly on $\dualJ$, the manifold where $X_h$ lives. However, one has to be cautious: it is well known that the Hamiltonian function $H$ on $\dualJ$ picks up extra terms under a time-dependent canonical transformation; in a way, there is no life for $X_h$ without the presence of $T^*E$, which therefore has to remain in the picture. Our aim is to explore in detail how conditions on $T^*E$ relate to objects on $\dualJ$ and vice versa. In Section~2, starting from a type $(1,1)$ tensor $R$ on $\dualJ$ with the property $R(dt)=0$, we recall the construction of the complete lifts of $R$ to both $T^*E$ and $\dualJ$ and their related properties. Of particular interest is that they both define a Poisson-Nijenhuis structure as soon as the Nijenhuis torsion of $R$ vanishes. Under an assumption of diagonalizability of $R$, the complete lift on $T^*E$ gives rise to an interesting distribution associated to any function $F$ on $T^*E$, and we characterize its integrability in Section~3. In the case of a function defining the image of a section $h:\dualJ\rightarrow T^*E$, there is a corresponding distribution on $\dualJ$. The interplay between the two distributions is studied in detail in Section~4. The integrability of both distributions is claimed to be an intrinsic version of Forbat's conditions for separability and we prove the claim in Section~5 by showing that we indeed recover Forbat's conditions in Darboux-Nijenhuis coordinates for the Poisson-Nijenhuis structures under consideration. An illustrative example is discussed in Section~6. Before starting, we should say that we owe a great deal of the inspiration for the distributions under consideration to a private meeting of one of us with Franco Magri back in 2001. To the best of our knowledge, Magri's ideas were never published, but they were the source of inspiration also for some of the results reported in \cite{Cra05} and \cite{Pedroni}.

\section{The complete lifts of a type $(1,1)$ tensor on $E$}

Let $R$ be a $(1,1)$ tensor on $E$ which vanishes on $dt$, in coordinates:
\begin{equation}
R = R^i_j(t,q)\fpd{}{q^i}\otimes dq^j + R^i_0(t,q)\fpd{}{q^i}\otimes dt. \label{R}
\end{equation}
It defines a fibre linear map on $T^*E$, given by
\begin{equation}
\tau_R: T^*E \rightarrow T^*E,\ (t,q^i,p_0,p_i)\mapsto (t,q^i,R^i_0p_i,R^i_jp_i). \label{Rmap}
\end{equation}
The complete lift of $R$ to $T^*E$, which we shall denote by $\cR$, is a well known construction. It can be defined for example by the relation (see \cite{CCS})
\begin{equation}
i_{\cR(Z)}d\theta_E = i_Z(\tau_R^*d\theta_E), \qquad \forall Z\in\vectorfields{T^*E}. \label{cR}
\end{equation}
In coordinates,
\begin{align}
\cR &= R^i_j \left(\fpd{}{q^i}\otimes dq^j + \fpd{}{p_j}\otimes dp_i\right) + R^i_0 \left(\fpd{}{q^i}\otimes dt + \fpd{}{p_0}\otimes dp_i\right) \nonumber\\
& \mbox{} + p_i\left(\fpd{R^i_j}{q^k} - \fpd{R^i_k}{q^j}\right) \fpd{}{p_j}\otimes dq^k  + p_i\left(\fpd{R^i_k}{t} - \fpd{R^i_0}{q^k}\right) \fpd{}{p_k}\otimes dt \nonumber \\
& \mbox{} + p_i\left(\fpd{R^i_0}{q^k} - \fpd{R^i_k}{t}\right)\fpd{}{p_0}\otimes dq^k. \label{cRc}
\end{align}
A direct construction of a complete lift of $R$ to $\dualJ$ is not that straightforward. In \cite{SW}, we developed a way to do it by an action on appropriately lifted vector fields. We will frequently refer to that paper for further properties and other lifting operations. For now, however, it suffices to introduce the complete lift $\wt{R}$ on $\dualJ$ by the property that it is $\rho$-related to $\cR$ on $T^*E$. That is to say, if $X\in\vectorfields{\dualJ}$ and $Y\in\vectorfields{T^*E}$ are any pair of $\rho$-related vector fields, then $\cR(Y)$ is $\rho$-related to $\wt{R}(X)$. In coordinates,
\begin{align}
\wt{R} &= R^i_j \left(\fpd{}{q^i}\otimes dq^j + \fpd{}{p_j}\otimes dp_i\right) + R^i_0 \fpd{}{q^i}\otimes dt \nonumber \\
& \mbox{} + p_i\left(\fpd{R^i_j}{q^k} - \fpd{R^i_k}{q^j}\right) \fpd{}{p_j}\otimes dq^k  + p_i\left(\fpd{R^i_k}{t} - \fpd{R^i_0}{q^k}\right) \fpd{}{p_k}\otimes dt, \label{tildeR}
\end{align}
from which the above cited link with $\cR$ is obvious.

If $N_R$ denotes the Nijenhuis torsion of $R$, which is a vector valued 2-form, we have the following important result.

\begin{thm}\  $N_{\cR}$ on $T^*E$ and $N_{\wt{R}}$ on $\dualJ$ vanish identically if and only if $N_R=0$.
\end{thm}
\proof\ The property is well known for the lift to a cotangent bundle (see \cite{YI} or \cite{CCS}). A proof for the lift to $\dualJ$ was given in \cite{SW}. \qed

A property of $\cR$ which immediately follows from the defining relation (\ref{cR}) is its symmetry with respect to $d\theta_E$, meaning that
\begin{equation}
d\theta_E(\cR(U),V) = d\theta_E(U,\cR(V)), \qquad \forall U,V\in\vectorfields{T^*E}. \label{sym}
\end{equation}
Obviously, $\cR$ then has a similar symmetry property with respect to the Poisson structure on $T^*E$ which is the inverse of $d\theta_E$. The manifold $\dualJ$ does not carry a symplectic structure, but it inherits a Poisson structure from $T^*E$ by projection \cite{GMS}. Let us denote the corresponding Poisson maps (i.e.\ the Poisson tensors regarded as map from 1-forms to vector fields) by $P_{T^{\!*}}$ and $P$ respectively.

\begin{thm}\ $(P_{T^{\!*}},\cR)$ and $(P,\wt{R})$ define Poisson-Nijenhuis structures on $T^*E$ and $\dualJ$ respectively, if and only if $N_R=0$.
\end{thm}
\proof\ There are basically three conditions for having a Poisson-Nijenhuis structure. One is the symmetry property referred to above (which equally holds for the projected Poisson structure on $\dualJ$). A second one is the vanishing of the so-called Magri-Morosi concommitant; this was proved to be the case in detail on $\dualJ$ in \cite{SW} and the proof for $T^*E$ can be found in \cite{CST}. The last condition is the vanishing of the Nijenhuis torsion, so that the statement then immediately follows from the preceding theorem. \qed

In what follows, we will assume firstly that the tensor $R$ is \emph{algebraically diagonalizable\/}, meaning that at each point $e\in E$, the endomorphism $R_e$ of $T_eE$ is diagonalizable, and secondly that \emph{the eigenvalues are distinct\/}. Since $R(dt)=0$, we know that zero is one of the eigenvalues, so that the remaining eigenvalues $\lambda_i$ (which generally will be functions of the coordinates $(t,q^i)$ of $e$) are non-zero by assumption. From the coordinate expressions (\ref{cRc}) and (\ref{tildeR}), it is clear that the coefficient matrices of $\cR$ and $\wt{R}$ have a block matrix structure. In particular, they have a $n\times n$ zero block corresponding to the lack of terms of the form $\partial/\partial q^i \otimes dp_j$, and have twice the block matrix $(R^i_j)$ along the diagonal. It then further readily follows that $\cR$ has double eigenvalues $(0,\lambda_i)$ and $\wt{R}$ has a single eigenvalue 0 and double eigenvalues $\lambda_i$. As a result, both matrices have
$\lambda \prod _{i=1}^n (\lambda-\lambda_i)$ as their minimal polynomial.

The main point about the extra assumption of diagonalizability is that, in the context of a Poisson-Nijenhuis structure, one can do better than merely algebraically diagonalize. Indeed, it is known that one can perform a diagonalization in coordinates which will at the same time produce Darboux coordinates for the Poisson tensor. In our situation, the natural bundle coordinates on $T^*E$ and $\dualJ$ already are Darboux coordinates for the Poisson structures under consideration. So we are prompted to be a little bit more careful about the construction of Darboux-Nijenhuis coordinates here. That is to say, we want to make sure that the diagonalization can be achieved by a time-dependent canonical transformation, i.e.\ the induced transformation of a time-dependent point transformation on $E$, which as such will not destroy the coordinate form of the Poisson tensor. We have proved in detail in \cite{SW} that this can be done for the Poisson-Nijenhuis structure $(P,\wt{R})$ on $\dualJ$, and it is an easy matter to verify that the same point transformation on $E$ will also induce a canonical transformation on $T^*E$ that does the job. Note further that in the course of proving the existence of an appropriate transformation $(t,q)\rightarrow (t,Q(t,q))$ on $E$, we found that each eigenvalue $\lambda_i$ will in the new coordinates at most depend on the corresponding coordinate $Q^i$. Hence we can formulate the following result.

\begin{thm}\ Let $R$ be a type $(1,1)$ tensor on $E$ which has the property $R(dt)=0$. Suppose that $N_R=0$ and that $R$ is algebraically diagonalizable with distinct eigenvalues. Then, there exists a coordinate transformation $(t,q)\rightarrow (t,Q(t,q))$ on $E$, which induces Darboux-Nijenhuis coordinates for both Poisson-Nijenhuis structures $(P_{T^{\!*}},\cR)$ and $(P,\wt{R})$ on $T^*E$ and $\dualJ$ respectively. In the new coordinates on $E$, $R$ takes the form \begin{equation}
R= \sum_{i=1}^n \lambda_{i}(Q^i) \fpd{}{Q^i} \otimes dQ^i. \label{diagR}
\end{equation}\label{DN}
\end{thm}
It follows that the coordinate expressions of $\cR$ and $\wt{R}$ in Darboux-Nijenhuis coordinates formally are identical and read,
\begin{equation}
\cR=\wt{R} =  \sum_{i=1}^n \lambda_{i}(Q^i) \left(\fpd{}{Q^i} \otimes dQ^i + \fpd{}{P_i}\otimes dP_i \right) . \label{diagcR}
\end{equation}

To conclude this section, we now show that the tensor $R$ also equips the manifold $\dualJ$ with a presymplectic structure. It will play a key role in the developments of the subsequent sections. In \cite{SW}, we defined a 1-form $\h{R}$ on $\dualJ$, called the \emph{horizontal lift of $R$\/}. Pointwise, its construction is determined by
\[
\langle X_m, \h{R}_m\rangle = \langle T\pi (X_m), R_{\pi(m)}(m)\rangle \quad \mbox{for all $m\in\dualJ,\ X\in\vectorfields{\dualJ}$}.
\]
[There is an unfortunate omission of the vector arguments in the definition formulated in \cite{SW}.]
In coordinates $\h{R}$ reads,
\begin{equation}
\h{R} = p_iR^i_j dq^j + p_iR^i_0 dt. \label{hR}
\end{equation}
Now consider the 2-form
\begin{align}
\omega_R := d\h{R} = h^*\tau^*_R d\theta_E&= R^j_i dp_j\wedge dq^i + R^j_0 dp_j\wedge dt \nonumber \\
& + \frac{1}{2} p_l\left(\fpd{R^l_j}{q^k} - \fpd{R^l_k}{q^j}\right) dq^k\wedge dq^j +
p_l\left(\fpd{R^l_0}{q^k} - \fpd{R^l_k}{t}\right) dq^k\wedge dt. \label{omegaR}
\end{align}
Clearly, $\omega_R$ is closed. In addition, the assumption about distinct eigenvalues of $R$ implies that $\det(R^i_j)\neq 0$. It is then clear from the above coordinate expression that $\omega_R$ has maximal rank, so that we have a presymplectic structure indeed. It is of some interest to have a look at the 1-dimensional kernel of $\omega_R$. If we compute $i_X\omega_R$ for an arbitrary $X\in\vectorfields{\dualJ}$, it is readily seen from the coefficients of $dp_j$ that $i_X\omega_R=0$ implies that $i_X \h{R}=0$. As a result $\Omega= \h{R}\wedge (\omega_R)^n$ cannot possibly be a volume form on $\dualJ$ so that $\h{R}$ does not define a contact structure.

\section{Distributions associated to functions on $T^*E$}

Under the assumptions of Theorem \ref{DN}, let $F$ be a function on $T^*E$. Consider the distribution
\[
\D{F} = \sp \,\{dF, \cR(dF), \cR^2(dF), \ldots, \cR^n(dF) \}^\circ
\]
consisting of the set of vector fields annihilating the indicated 1-forms. The sequence of these 1-forms certainly breaks down at the power $n$ of $\cR$, because the minimal polynomial of $\cR$ has degree $n+1$. Assume that $F$ and $R$ are such that the defining 1-forms are linearly independent (except for isolated points), so that $\D{F}$ also has dimension $n+1$ and will be Lagrangian, provided it is isotropic or co-isotropic.

\begin{lemma}\ The orthogonal complement $\Dperp{F}$ of\ $\D{F}$ is given by
\[
\Dperp{F} = \sp\,\{X_F,\cR(X_F), \ldots, \cR^n(X_F)\}.
\]
\end{lemma}
\proof\ For all $Y\in\D{F}$ and $k=0,\ldots, n$ we have, using the symmetry of $\cR$ with respect to $\omega_E$, that
\[
\omega_E(\cR^k(X_F),Y) = \omega_E(X_F,\cR^k(Y)) = - dF(\cR^k(Y)) = - \cR^k(dF)(Y)=0,
\]
which shows that $\cR^k(X_F)$ belongs to $\Dperp{F}$ for $k=0,\ldots,n$. Moreover, in open domains where the defining 1-forms of $\D{F}$ are linearly independent, the same is true for the vector fields $\cR^k(X_F)$. Indeed, if $\sum_{k=0}^n a_k \cR^k(X_F)=0$,
we have for all $Z\in\vectorfields{T^*E}$,
\[
0= \omega_E\Big(\sum_{k=0}^n a_k \cR^k(X_F),Z\Big) = - \sum_{k=0}^n a_k \langle \cR^k(Z), dF\rangle = - \Big\langle Z, \sum_{k=0}^n a_k \cR^k(dF)\Big\rangle
\]
which implies that all functions $a_k$ must be zero. By dimension, therefore, the $\cR^k(X_F)$ span $\Dperp{F}$. \qed

\begin{lemma} \ The distribution $\D{F}$ is Lagrangian, i.e.\ $\D{F}=\Dperp{F}$.
\end{lemma}
\proof\ We will show that $\D{F}$ is co-isotropic. For that purpose, consider
\[
\Big\langle \cR^l(X_F),\cR^k(dF)\Big\rangle = \Big\langle \cR^{k+l}(X_F),dF\Big\rangle = - \omega_E\big(X_F, \cR^{k+l}(X_F)\big).
\]
Again, by the symmetry of $\cR$ with respect to $\omega_E$, we have
\[
\omega_E\big(X_F, \cR^{k+l}(X_F)\big) = \omega_E\big(\cR^{k+l}(X_F), X_F\big),
\]
but then the skew-symmetry of $\omega_E$ implies that this is identically zero. Since this is valid for all $l$ and $k$, we conclude from the first line that $\Dperp{F}\subseteq D_F$. The dimension then implies that we have equality and thus a Lagrangian distribution. \qed

In what follows, we will denote the distribution simply by $\D{F}$ even when we appeal to the defining relation of $\Dperp{F}$. Note further that it follows from both defining relations and the degree of the minimal polynomial of $\cR$ that $\cR(\D{F})\subset \D{F}$.

Naturally, we are interested in the case that $\D{F}$ is Frobenius integrable. In preparation of our main theorem about this integrability, we list the following general property of derivations.

\begin{lemma}\
Let $L$ be a type $(1,1)$ tensor field on an arbitrary manifold $M$. Then, for any 1-form $\alpha$ and vector fields $X,Y$ on $M$:
\begin{equation}
d_L(L\alpha)(X,Y) = d\alpha (LX,LY) + \alpha\big(N_L(X,Y)\big). \label{lem3}
\end{equation}
\end{lemma}
\proof\ Since $d_L=i_Ld-di_L$, we have
\[
d_L(L\alpha)(X,Y) = d(L\alpha)(LX,Y) + d(L\alpha)(X,LY) - d(L^2\alpha)(X,Y).
\]
Using the general property $d\alpha(X,Y)=\lie{X}(\alpha(Y))-\lie{Y}(\alpha(X)) - \alpha([X,Y])$, this easily reduces to
\begin{align*}
d_L(L\alpha)(X,Y) &= \lie{LX}((L\alpha)(Y)) - \lie{LY}((L\alpha)(X)) \\
& \qquad - (L\alpha)\big([LX,Y]\big) - (L\alpha)\big([X,LY]\big) + (L^2\alpha)\big([X,Y]\big), \\
&= \lie{LX}(\alpha(LY)) - \lie{LY}(\alpha(LX)) - \alpha\big([LX,LY]\big) + \alpha\big(N_L(X,Y)\big),
\end{align*}
from which the result now follows. \qed

In particular, if $L$ has vanishing Nijenhuis torsion, then
\[
d_L(L\alpha)(X,Y) = d\alpha (LX,LY).
\]

\begin{thm}
Let $R$ be a $(1,1)$ tensor on $E$ with the property $R(dt)=0$, which is algebraically diagonalizable with distinct eigenvalues and has vanishing Nijenhuis torsion. Suppose $F$ is a function on $T^*E$ for which the defining 1-forms of the distribution $\D{F}$ are linearly independent. Then $\D{F}$ is integrable if and only if $d d_{\cR}F\big|_{\D{F}} =0$.
\label{integrableD}
\end{thm}
\proof\ Looking at the defining co-distribution of $\D{F}$ and putting $\alpha_i = d_{\cR^i}F$ for shorthand, we know that $\D{F}$ is integrable if and only if $d\alpha_i = \sum_{l=0}^n \theta^l_i\wedge \alpha_l$ for some 1-forms $\theta^l_i$. By extending the $\alpha_i$ to a local basis for $\oneforms{T^*E}$, it is easy to see that this is further equivalent to $d\alpha_i\big|_{\D{F}}=0$ for all $i$. Hence, if $\D{F}$ is integrable, we have in particular that $d\alpha_1\big|_{\D{F}}= d d_{\cR}F\big|_{\D{F}} =0$.

Conversely, assuming $d d_{\cR}F\big|_{\D{F}} =0$, we first observe that $d\alpha_0\big|_{\D{F}}=0$ since $d\alpha_0=0$, and also
\[
d_{\cR}\alpha_0\big|_{\D{F}} = - d i_{\cR}\alpha_0\big|_{\D{F}} = - d d_{\cR}F\big|_{\D{F}} =0.
\]
We now proceed further by induction. Assuming that $d\alpha_i\big|_{\D{F}}=0$ and $d_{\cR}\alpha_i\big|_{\D{F}} = 0$, we will show that the same properties hold for $\alpha_{i+1}$. Firstly, for all $X,Y\in\D{F}$, using (\ref{lem3}) and the fact that $N_{\cR}=0$, we conclude that
\[
d_{\cR}\alpha_{i+1}(X,Y) = d_{\cR}(\cR\alpha_i)(X,Y) = d\alpha_i(\cR X,\cR Y)=0,
\]
since $\cR(\D{F})\subset \D{F}$. Secondly,
\[
d\alpha_{i+1}(X,Y) = d(\cR\alpha_i)(X,Y) = i_{\cR}d\alpha_i(X,Y) - d_{\cR}\alpha_i(X,Y),
\]
which reduces to the first term on the right by the induction hypothesis and then in fact to zero in view of $\cR(\D{F})\subset \D{F}$ and the induction hypothesis again. The conclusion is that, in particular, $d\alpha_i\big|_{\D{F}}=0$ for all $i$ and hence that $\D{F}$ is integrable. \qed

There is a direct link, which we will briefly sketch now, between the integrability of $\D{F}$ and the classical Hamilton-Jacobi equation for the Hamiltonian $F$ on $T^*E$. It is well known that in the neighbourhood of a regular point, i.e.\ a point where the distribution $\D{F}$ is transversal to the fibers, every Lagrangian submanifold of $T^*E$ is the image of the differential of a function defined on an open subset of $E$ (see e.g.\ \cite{BT} or \cite{LibMar}, Appendix~7). In such a regular point, the vectors spanning $\D{F}^\perp$ will be linearly independent if and only if their projections onto $E$ are linearly independent. If $\pi_E= \pi\circ\rho$ denotes the projection of $T^*E$ onto $E$, we have at each point $(t,q,p_0,p)$ of $T^*E$ that
\begin{align*}
T\pi_E(X_{F}(t,q,p_0,p)) &= \fpd{F}{p_0}\left.\fpd{}{t}\right|_{(t,q)} + \fpd{F}{p_i}\left.\fpd{}{q^i}\right|_{(t,q)} \\
T\pi_E\big(\cR(X_{F})(t,q,p_0,p)\big) &= \Big(R^i_j\fpd{F}{p_j} + R^i_0\fpd{F}{p_0}\Big)\left.\fpd{}{q^i}\right|_{(t,q)} \\
&\vdots \\
T\pi_E\big(\cR^n(X_{F})(t,q,p_0,p)\big) &= {R^{(n-1)}}^i_l \Big(R^l_j\fpd{F}{p_j} + R^l_0\fpd{F}{p_0}\Big)\left.\fpd{}{q^i}\right|_{(t,q)}.
\end{align*}
If we think of $T\pi_E(X_{F}(t,q,p_0,p))$ as being expressed in terms of the basis of eigenvectors of $R$, it is clear that independence of the above set requires that none of the coefficients in that expression be zero.
Equivalently, this means that $\partial{F}/\partial p_0 \neq 0$ and the vector with components $R^i_j(\partial F/\partial p_j) + R^i_0 (\partial F/\partial p_0)$ is spanned by all eigenvectors of the matrix $(R^i_j)$.
If $\D{F}$ is integrable, we have a foliation of $T^*E$ into Lagrangian submanifolds and in the neighbourhood of regular points each of them will be generated by the differential of a function $S$, which therefore locally extends to a generating function defined on $T^*E$. Moreover, since $dF|_{\D{F}}=0$ by construction, we further conclude that on such a submanifold, $S$ will be a solution of the partial differential equation
\[
F\Big(t,q^i,\fpd{S}{t},\fpd{S}{q^i}\Big) = \mbox{constant}
\]
which is the Hamilton-Jacobi equation for $F$.
For another excellent and extensive account of the geometry of the Hamilton-Jacobi equation we refer to \cite{MMM}. Observe though that in all references cited in this context, the base manifold $E$ for the time-dependent case is taken to be a product manifold $Q\times\R$.

This brings us to the point that the above type of Hamilton-Jacobi equation is of course not exactly what we are interested in: it is in some sense Hamilton-Jacobi theory for an autonomous Hamiltonian where one of the position coordinates happens to be denoted by $t$. The case of interest is when $F$ is of the form $F=\wt{H}:= p_0 + H(t,q^i,p_i)$, and hence $\wt{H}=0$ defines a section of $\rho:T^*E\rightarrow \dualJ$ and a time-dependent Hamiltonian system on $\dualJ$. The corresponding Hamilton-Jacobi equation then is of the form:
\begin{equation}
\fpd{S}{t} + H\Big(t,q^i,\fpd{S}{q^i}\Big) = 0. \label{HJ}
\end{equation}
For this case, we study the existence of a corresponding integrable distribution on $\dualJ$ in the next section.

\section{A corresponding distribution on $\dualJ$ for given section $h:\dualJ\rightarrow T^*E$}

Consider a section $h:\dualJ\rightarrow T^*E$, whose image is the set of points in $T^*E$ such that $\wt{H}=0$. We have a corresponding Hamiltonian vector field $X_{\wt{H}}$ on $T^*E$, with coordinate expression,
\begin{equation}
X_{\wt{H}} = \fpd{}{t} + \fpd{H}{p_i}\fpd{}{q^i} - \fpd{H}{q^i}\fpd{}{p_i} - \fpd{H}{t}\fpd{}{p_0}, \label{XH}
\end{equation}
and the Hamiltonian vector field $X_h$ on $\dualJ$ as in (\ref{Xh}). Clearly, $X_{\wt{H}}$ projects onto $X_h$, in other words $X_{\wt{H}}$ and $X_h$ are $\rho$-related. We have already mentioned in Section~2 that also the tensor fields $\cR$ and $\wt{R}$ are $\rho$-related. As a result, if we consider the distribution $\D{h}$ on $\dualJ$, defined by
\begin{equation}
\D{h} = \sp\, \{X_h, \wt{R}(X_h), \wt{R}^2(X_h),\ldots, \wt{R}^n(X_h)\}, \label{Dh}
\end{equation}
it is clear that $\D{\wt{H}}$ on $T^*E$ and $\D{h}$ on $\dualJ$ are $\rho$-related. We wish to show now that, more importantly, they are also $h$-related. A vector field $X$ on $\dualJ$ is $h$-related to $Y$ on $T^*E$ if $Th\circ X = Y\circ h$ or equivalently $Y(F)\circ h = X(F\circ h)$ for all functions $F$ on $T^*E$. It is easy to verify that this translates into the following relationship between the local coordinate expressions. If
\[
X = X^0 \fpd{}{t} + X^i\fpd{}{q^i} + U_i\fpd{}{p_i},
\]
then, with a little abuse of notation because we omit the pullback under $\rho$ for functions which come from $\dualJ$, $Y$ will necessarily be of the form
\begin{equation}
Y = X^0 \fpd{}{t} + X^i\fpd{}{q^i} + U_i\fpd{}{p_i} - X(H)\fpd{}{p_0}= X - X(H)\fpd{}{p_0}. \label{hrel}
\end{equation}
\begin{lemma}\ $X\in\vectorfields{\dualJ}$ and $Y\in\vectorfields{T^*E}$ are $h$-related if and only if $Y$ projects onto $X$ and $Y(\wt{H})=0$.
\end{lemma}
\proof\ This is immediately clear from the above coordinate expressions. More intrinsically, $Y$ must project onto $X$ because $\rho$ and $h$ are each others inverse when restricted to the image of $h$ and $Y(\wt{H})=0$ reflects the fact that $Y$ must be tangent to this image. \qed

Obviously, $X_{\wt{H}}(\wt{H})=0$, hence $X_h$ and $X_{\wt{H}}$ are $h$-related. But it is certainly not true that also the tensor fields $\wt{R}$ and $\cR$ would be $h$-related, i.e.\ that they would map general $h$-related vector fields into $h$-related vector fields. That is true, however, when we restrict ourselves to the sequence of vector fields defining the distributions $\D{h}$ and $\D{\wt{H}}$.

\begin{lemma}\ The vector fields $\wt{R}^k(X_h)\in\vectorfields{\dualJ}$ and $\cR^k(X_{\wt{H}})\in\vectorfields{T^*E}$ are $h$-related for all $k$.
\end{lemma}
\proof\ We already know that the vector fields under consideration are $\rho$-related. In addition, by the fact that $\D{\wt{H}}= \D{\wt{H}}^\perp$, we have: $\langle \cR^k(X_{\wt{H}}),d\wt{H}\rangle=0$ for all $k$.
\qed

As before, we assume that $\wt{H}$ and $R$ are such that $\D{\wt{H}}$ has dimension $n+1$ in open domains of $T^*E$. It is clear then, by a pointwise projection argument for example, that the defining vector fields of $\D{h}$ in (\ref{Dh}) are also linearly independent, so that $\D{h}$ is a distribution of dimension $n+1$ on $\dualJ$.

\begin{thm}\ $\D{h}$ is an integrable distribution on $\dualJ$ if and only if $\D{\wt{H}}$ is integrable on $T^*E$.
\end{thm}
\proof\ If two vector fields on $\dualJ$ are $h$-related to corresponding vector fields on $T^*E$, then so are their Lie brackets. By way of example, consider the pair $(X_h,\wt{R}(X_h))$ on $\dualJ$ and the corresponding pair $(X_{\wt{H}},\cR(X_{\wt{H}}))$ on $T^*E$, but the reasoning below applies just as well to any other pair. The fact that their brackets are also $h$-related means that, in terms of the simplified notations used in (\ref{hrel}), we have
\[
\big[X_{\wt{H}},\cR(X_{\wt{H}})\big] = \big[X_h,\wt{R}(X_h)\big] - \big[X_h,\wt{R}(X_h)\big](H)\fpd{}{p_0}.
\]
Now, if $\D{h}$ is integrable, we have that $\big[X_h,\wt{R}(X_h)\big]= \sum_{k=0}^n a_k \wt{R}^k(X_h)$ for some functions $a_k$ on $\dualJ$. Using this in the above equality, the right-hand side clearly becomes, again by the formal general rule (\ref{hrel}), the expression for the $h$-related vector field $\sum_{k=0}^n a_k \cR^k(X_{\wt{H}})$. This shows that the bracket $\big[X_{\wt{H}},\cR(X_{\wt{H}})\big]$ belongs to $\D{\wt{H}}$, and similarly for all other pairs, so that $\D{\wt{H}}$ is integrable. Conversely, assume that $\D{\wt{H}}$ is integrable, then $\big[X_{\wt{H}},\cR(X_{\wt{H}})\big] = \sum_{k=0}^n a_k \cR^k(X_{\wt{H}})$, for some $a_k$ which in principle could be functions on $T^*E$. But all vector fields $\cR^k(X_{\wt{H}})$ in that sum are $h$-related to a corresponding element of $\D{h}$, so that
\[
\sum_{k=0}^n a_k \cR^k(X_{\wt{H}}) = \sum_{k=0}^n a_k \wt{R}^k(X_h) - \sum_{k=0}^n a_k \wt{R}^k(X_h)(H)\fpd{}{p_0}.
\]
Identifying the right-hand sides of both displayed equalities, we conclude that necessarily
$\big[X_h,\wt{R}(X_h)\big]= \sum_{k=0}^n a_k \wt{R}^k(X_h)$. The left-hand side in this relation manifestly is a vector field on $\dualJ$, so that there cannot be any $p_0$-dependence in the overall expression on the right. Therefore, if some of the $a_k$ would explicitly depend on $p_0$, the partial sum of such terms on the right would have to vanish. But if vector fields such as the $\wt{R}^k(X_h)$ are linearly independent as vector fields on $\dualJ$, then they are also linearly independent as vector fields along the projection $\rho:T^*E\rightarrow \dualJ$. This implies that all $a_k$ in that partial sum eventually must vanish. So in the end, we will have an equality of the form $\big[X_h,\wt{R}(X_h)\big]= \sum_{k=0}^n a_k \wt{R}^k(X_h)$ with $a_k$ which, without loss of generality, can be seen as functions on $\dualJ$. Repeating this argument for all possible brackets of vector fields of the form $\cR^k(X_{\wt{H}})$ will lead us to the conclusion that also $\D{h}$ is integrable.  \qed

By Theorem \ref{integrableD}, integrability of $\D{\wt{H}}$ is reduced to the condition $dd_{\cR}{\wt{H}}\big|_{\D{\wt{H}}} =0$, and we now know that this will equally ensure integrability of $\D{h}$. But there is no doubt that it would be more satisfactory still to characterize integrability of $\D{h}$ by a condition expressed in terms of objects living on $\dualJ$. This is our final goal for this section and it will be achieved with the aid of the presymplectic structure on $\dualJ$ defined by $\omega_R$.

\begin{lemma}
The presymplectic form $\omega_R$, defined by $\omega_R= d\h{R}$ has the additional property that
\begin{equation}
\omega_R = h^*\tau^*_R d\theta_E. \label{omegaR2}
\end{equation}
\end{lemma}
\proof\ It suffices to verify from (\ref{Rmap}) and (\ref{hR}) that $h^*\tau^*_R \theta_E = \h{R}$ indeed. \qed

Since $\tau^*_R d\theta_E$ is the 2-form needed to define $\cR$ (see (\ref{cR})), the idea now is to transfer certain properties from $T^*E$ to $\dualJ$ by pulling back via $h$. Of course, such a pullback works well for forms, but is in general not well defined when it concerns the contraction of a form with an arbitrary vector field. But it does work when the vector fields involved have an $h$-related companion on $\dualJ$, as we briefly recall first in a general setting.

Let $h$ be a smooth map from a manifold $M$ into a manifold $N$ and let $Y\in\vectorfields{N}$ be $h$-related to $X\in\vectorfields{M}$, so that $Th\circ X= Y\circ h$. Then, for any form $\omega\in\forms{k}{N}$, we can define $h^*(i_Y\omega)\in\forms{k-1}{M}$ as follows. For any $m\in M$ and $v_1,\ldots, v_{k-1}\in T_mM$, put
\begin{align*}
h^*(i_Y\omega)(m)(v_1,\ldots,v_{k-1}) &= (i_Y\omega)(h(m))(Th(v_1),\ldots,Th(v_{k-1})) \\
&= \omega(h(m))(Th(X_m),Th(v_1),\ldots,Th(v_{k-1})) \\
&= (h^*\omega)(m)(X_m, v_1,\ldots,v_{k-1}),
\end{align*}
from which it follows that $h^*(i_Y\omega) = i_X(h^*\omega)$.

In the course of the proof of Lemma~1, applied to the case of $\D{\wt{H}}$, we have seen that
\begin{equation}
i_{\cR^k(X_{\wt{H}})}d\theta_E = - \cR^k(d\wt{H}), \qquad \mbox{for all\ } k. \label{aux1}
\end{equation}
On the other hand, we know from the defining relation (\ref{cR}) of $\cR$ that
\begin{equation}
i_{\cR^k(X_{\wt{H}})}d\theta_E = i_{\cR^{k-1}(X_{\wt{H}})}\tau^*_Rd\theta_E. \label{aux2}
\end{equation}
It follows that
\begin{equation}
i_{\cR^{k-1}(X_{\wt{H}})}\tau^*_Rd\theta_E = - \cR^k(d\wt{H}), \qquad k=1,\ldots, n. \label{aux3}
\end{equation}

\begin{thm}\ The distribution $\D{h}$ on $\dualJ$ is integrable if and only if
\begin{equation}
\left.\lie{X_h}\omega_R \right|_{\D{h}} = 0. \label{intDh}
\end{equation}
\end{thm}
\proof\ By Lemma~5 about $h$-related vector fields, we can pull back the relations (\ref{aux3}) under $h$, to obtain
\[
i_{\wt{R}^{k-1}(X_h)}\omega_R = - h^*\Big(\cR^k(d\wt{H})\Big), \qquad k=1,\ldots, n.
\]
Taking the exterior derivative of this relation for the case $k=1$ and knowing that $\omega_R$ is closed, the result now immediately follows from Theorem \ref{integrableD}, applied to the case where $F=\wt{H}$. This final step of course again relies on the fact that we have bases of $\D{\wt{H}}$ and $\D{h}$ consisting of $h$-related vector fields. \qed

\section{Darboux-Nijenhuis coordinates and Forbat's conditions for separability of the Hamilton-Jacobi equation}

As announced in the introduction, the integrability of $\D{h}$ on $\dualJ$ (or equivalently of $\D{\wt{H}}$ on $T^*E$) is claimed to be an intrinsic formulation of Forbat's conditions for separability of the time-dependent Hamilton-Jacobi equation. So, we should be able to show that there exists a selection of natural coordinates, such that the condition (\ref{intDh}) on $\dualJ$ (or equivalently the condition $d d_{\cR}\wt{H}\big|_{\D{\wt{H}}} =0$ on $T^*E$ coming from Theorem \ref{integrableD}) precisely reproduces the conditions (\ref{Forbat1},\ref{Forbat2}). Needless to say, if such a preferred coordinate system exists, it should have made its appearance in the course of the theoretical developments. It should therefore not be a surprise that we actually claim that Darboux-Nijenhuis coordinates on $\dualJ$ or $T^*E$ do the job. We shall show this for the condition (\ref{intDh}), but it can equally well be carried out for the equivalent condition on $T^*E$.

Suppose we have found the coordinate transformation $(t,q)\rightarrow (t,Q(t,q))$ which diagonalizes the tensor $R$ on $E$, and let $(t,q,p)\rightarrow (t,Q(t,q),P(t,q,p))$ be the induced time-dependent canonical transformation on $\dualJ$. In other words, we have that
\[
P_i(t,q,p) = p_l\fpd{q^l}{Q^i}(t,Q(t,q)).
\]
Then, we know by Theorem \ref{DN} that $\wt{R}$ will take the form (\ref{diagcR}). At the same time, the Hamiltonian vector field $X_h$ will have changed its appearance: explicitly, if $H(t,q,p)$ was the Hamiltonian function in the original coordinates, the Hamiltonian $K(t,Q,P)$ in the new coordinates will be given by (from the induced transformation of $p_0$ on $T^*E$)
\begin{equation}
K(t,Q,P) = H - p_l\fpd{q^l}{t}, \label{K}
\end{equation}
with the understanding that the right-hand side has to be expressed in terms of the new variables. If we now compute the vector fields $\wt{R}^k(X_h)$ spanning the distribution $\D{h}$, we readily observe that
\begin{align*}
  \left(\begin{array}{c}X_{h}\vspace{0.15cm}\\ \wt{R}(X_h)\vspace{0.15cm}\\ \wt{R}^2(X_h)\vspace{0.15cm}\\ \vdots\vspace{0.15cm}\\ \wt{R}^n(X_h)\end{array}\right)
  &=\left(\begin{array}{ccccc}
  1&1&1&\ldots&1\vspace{0.15cm}\\ 0&\lambda_1&\lambda_2&\ldots&\lambda_n\vspace{0.15cm}\\
  0&\lambda_1^2&\lambda_2^2&\ldots&\lambda_n^2\vspace{0.15cm}\\ \vdots&\vdots&\vdots&\ddots&\vdots\vspace{0.15cm}\\ 0&\lambda_1^n&\lambda_2^n&\ldots&\lambda_n^n\end{array}\right)
  \left(\begin{array}{c}\fpd{}{t}\vspace{0.15cm}\\ \fpd{K}{P_1}\fpd{}{Q^1}-\fpd{K}{Q^1}\fpd{}{P_1}\vspace{0.15cm}\\  \fpd{K}{P_2}\fpd{}{Q^2}-\fpd{K}{Q^2}\fpd{}{P_2}\vspace{0.15cm}\\ \vdots\vspace{0.15cm}\\ \fpd{K}{P_n}\fpd{}{Q^n}-\fpd{K}{Q^n}\fpd{}{P_n}\end{array}\right).
\end{align*}
This strongly suggests a change of basis, which should no doubt simplify the calculations for the condition (\ref{intDh}). So we put (no sum)
\begin{equation}
X_0=\fpd{}{t},\quad X_i=\fpd{K}{P_i}\fpd{}{Q^i}-\fpd{K}{Q^i}\fpd{}{P_i}, \;\; i=1,\ldots,n \label{eigenX}
\end{equation}
and observe that this is in fact a set of eigenvectors for the tensor $\wt{R}$, as given by (\ref{diagcR}) in the coordinates under consideration. There is more to say about this observation. Since $\cR$ formally is identical to $\wt{R}$ in Darboux-Nijenhuis coordinates, a similar computation of the vector fields $\cR^k(X_{\wt{H}})$ which span the distribution $\D{\wt{H}}$ on $T^*E$, will generate via the same non-singular transition matrix a basis of eigenvectors for $\cR$, given by (no sum)
\begin{equation}
Y_0=\fpd{}{t}-\fpd{K}{t}\fpd{}{P_0},\quad Y_i=\fpd{K}{P_i}\fpd{}{Q^i}-\fpd{K}{Q^i}\fpd{}{P_i}, \;\; i=1,\ldots,n. \label{eigenY}
\end{equation}
It is further interesting to notice that the vector fields $X_k$ on $\dualJ$ and $Y_k$ on $T^*E$ are $h$-related. Hence, we have proved, by passing to a special selection of coordinates, the following useful addition to Lemma~5.

\begin{lemma}\ There exists a transition to local bases for $\D{h}$ and $\D{\wt{H}}$ which consist of eigenvectors of $\wt{R}$ and $\cR$ respectively, and preserves the property that the generating vector fields are $h$-related. \qed
\end{lemma}

Let us now finally express the condition (\ref{intDh}) for integrability of $\D{h}$ by making use of the new basis of eigenvectors $X_k$ ($k=0,\ldots,n$). From (\ref{omegaR}) and (\ref{diagR}), we see that in the new coordinates, $\omega_R$ takes the simple form
\[
\omega_R= \sum_{l=1}^n \lambda_l(Q^l) dP_l\wedge dQ^l,
\]
from which it follows that
\begin{align*}
\lie{X_h}\omega_R &= \sum_{l=1}^n \fpd{K}{P_l} \frac{d\lambda_l}{dQ^l} dP_l\wedge dQ^l + \sum_{l=1}^n \lambda_l\left(\fpdd{K}{Q^l}{t} dQ_l\wedge dt + \fpdd{K}{P_l}{t} dP_l\wedge dt\right) \\
& \hspace*{-1cm} + \sum_{k,l} (\lambda_l-\lambda_k) \left(\fpdd{K}{P_l}{Q^k} dP_l\wedge dQ^k + \frac{1}{2}\Big(\fpdd{K}{Q^k}{Q^l}dQ^l\wedge dQ^k + \fpdd{K}{P_k}{P_l}dP_l\wedge dP_k\Big)\right).
\end{align*}
The first term does not contribute anything when acting on the basis of eigenvectors (\ref{eigenX}). From the second term, it follows that
\[
\lie{X_h}\omega_R(X_i,X_0) = \lambda_i \left(\fpdd{K}{Q^i}{t} \fpd{K}{P_i}- \fpdd{K}{P_i}{t}\fpd{K}{Q^i}\right)\qquad \mbox{(no sum)}.
\]
The last term implies that
\begin{align*}
\lie{X_h}\omega_R(X_i,X_j) &= (\lambda_j-\lambda_i)\left(\fpdd{K}{P_i}{Q^j}\fpd{K}{Q^i}\fpd{K}{P_j} + \fpdd{K}{P_j}{Q^i}\fpd{K}{Q^j}\fpd{K}{P_i} \right. \\
& \quad \left. - \fpdd{K}{Q^i}{Q^j} \fpd{K}{P_i}\fpd{K}{P_j} -
\fpdd{K}{P_i}{P_j} \fpd{K}{Q^i}\fpd{K}{Q^j} \right) \qquad \mbox{(no sums)}.
\end{align*}
Since the $\lambda_i$ are nonzero and distinct, it is clear now that $\left.\lie{X_h}\omega_R \right|_{\D{h}}=0$ precisely gives rise to the Forbat conditions (\ref{Forbat1},\ref{Forbat2}).

We summarize our main results in the following theorem.

\begin{thm} Let $E$ be a bundle over $\R$ of dimension $n+1$. Let $h$ be a section of the bundle $\rho:T^*E\rightarrow \dualJ$ and denote by $X_h$ the corresponding Hamiltonian vector field on $\dualJ$. Let $R$ be a type $(1,1)$ tensor field on $E$ with the following properties: (i) $R(dt)=0$, (ii) $N_R=0$, (iii) $R$ is algebraically diagonalizable with distinct eigenvalues. Assume further that the $n+1$ vector fields $X_h, \wt{R}(X_h), \ldots \wt{R}^n(X_h)$ are linearly independent, where $\wt{R}$ is the complete lift of $R$ to $\dualJ$. Consider the 2-form $\omega_R=d\h{R}$, where $\h{R}$ is the horizontal lift of $R$ to $\dualJ$. Then the distribution $\D{h}= \sp\,\{X_h, \wt{R}(X_h), \ldots \wt{R}^n(X_h)\}$ is Lagrangian, and is integrable provided that $\left.\lie{X_h}\omega_R \right|_{\D{h}}=0$. In Darboux-Nijenhuis coordinates for the Poisson-Nijenhuis structure which $\wt{R}$ defines on $\dualJ$, these integrability conditions are exactly Forbat's necessary and sufficient conditions for separability of a time-dependent Hamilton-Jacobi equation.
\end{thm}

\section{An example}

Having obtained an intrinsic characterization of Hamilton-Jacobi separability in the form of a set of conditions which in principle can be tested in any coordinate system, i.e.\ prior to knowing separation coordinates, the next challenging question is of course: ``What is the practical content of these conditions?". In other words, if one comes along with a given (time-dependent) Hamiltonian, how should we proceed to test whether a change of coordinates exists which will transform the Hamiltonian into one which is Hamilton-Jacobi separable? A thorough investigation of this problem is work for the future: we believe that it should be possible to make progress concerning an algorithmic path for testing existence of separation coordinates and constructing them, by making use of the classification results available in reference \cite{Frans} for example. For the sake of presenting an illustrative example, however, we can limit ourselves to a less ambitious goal.

There are examples in the literature of case studies, mostly in the context of autonomous systems though, where the starting point is a Hamiltonian system with a known separable (or, more generally, integrable) potential, and the purpose is to find a broader class of separable potentials. It is then customary to make a certain ansatz about the form of the potentials one is looking for. For our example, we start from a time-dependent Hamiltonian (with $n=2$) of the form
\begin{align}
H &= \onehalf p_1^2 + \onehalf t\,p_2^2 + 2t^3(t\,q_1^2 + q_2^2) + c_1(t)p_1q_1 + c_2(t)p_2q_2 \nonumber \\[1mm]
& \mbox{} \qquad + a_1(t)q_1^3 + a_2(t) q_1^2q_2 + a_3(t) q_1q_2^2 + a_4(t)q_2^3. \label{exH}
\end{align}
It contains six as yet arbitrary functions of time and the point about the explicitly specified terms is that they form part of a Hamiltonian with quadratic potential which we know to be separable after a linear change of coordinates. [For clarity, we use lower indices for the $q$-variables in such explicit polynomial expressions.] We shall likewise make an ansatz concerning the tensor $R$ whose complete lift should help to identify separable Hamiltonians through the conditions expressed in our main theorem. We take $R$ to be of the form
\begin{align}
R &= t q_1 \left(\fpd{}{q_1}\otimes dq_1 + \fpd{}{q_2}\otimes dq_2\right) + q_2\left(\sigma_1(t)\fpd{}{q_1}\otimes dq_2 + \sigma_2(t)\fpd{}{q_2}\otimes dq_1\right) \nonumber \\
& \mbox{} \qquad + R^1_0(t,q)\fpd{}{q_1}\otimes dt +  R^2_0(t,q)\fpd{}{q_2}\otimes dt. \label{exR}
\end{align}
Again, the first term is completely specified and is inspired by the quadratic case we know, but we still have four arbitrary functions at our disposal. In fact, it turns out that one could allow for an arbitrary function of time in the first term as well, but it would inevitably be forced to equal $t$ later in the process. Concerning the specified $q$-dependence in the leading coefficients $R^i_j(t,q)$ in (\ref{exR}), the motivation is the following. We insist on staying within the category of Hamiltonians which become separable after a linear change of coordinates. The feature of the tensor $R$ in (\ref{exR}) which will guarantee this is that the eigenfunctions of the matrix $R^i_j$ are given by $\lambda_i = tq_1 \pm \sqrt{\sigma_1\sigma_2}\,q_2$. As we learn from the general expression (\ref{diagR}), the coefficients of a suitable $R$ in Darboux-Nijenhuis coordinates can be arbitrary functions $\lambda_i(Q^i)$, each depending on a single coordinate $Q^i$. Hence, they can be chosen to be the $Q^i$ themselves. The structure of the $R^i_j$ in (\ref{exR}) therefore already tells us what the linear change of coordinates will be after pinning down the freedom in (\ref{exH}) and (\ref{exR}). Note further that we have to require that none of the $\sigma_i(t)$ becomes zero, since we want $R$ to have distinct eigenvalues. The ensuing coordinate change will of course be valid in the domain where $\sigma_1\sigma_2$ is positive.

The first condition we impose now is that $R$ should have vanishing Nijenhuis torsion. This turns out to fix the function $\sigma_2$ as
\[
\sigma_2(t) =t.
\]
The conditions which involve the as yet arbitrary $R^i_0$-components read (with partial derivatives denoted by a comma, and summation over $k$):
\[
R^i_k(R^k_{0,j} - R^k_{j,t}) = R^k_jR^i_{0,k} - R^k_0R^i_{j,k}.
\]
In our case, this is a set of four partial differential equations for the $R^i_0$, but quite remarkably, they combine in such a way that the $R^i_0$ can be determined from algebraic relations. This is in fact an interesting feature which is a consequence of our restriction to tensors which satisfy $R(dt)=0$. Taking $\sigma_2=t$ and $\sigma_1\neq 0$ into account, we obtain that
\begin{align*}
R^1_0 &=  q_1^2 + \onehalf \big((\sigma_1/t) + \dot{\sigma}_1\big) q_2^2, \\
R^2_0 &= \onehalf \big( (t\dot{\sigma}_1/\sigma_1) + 3\big) q_1q_2.
\end{align*}
We can now compute the vector fields $X_h,\wt{R}(X_h),\wt{R}^2(X_h)$ spanning the distribution $\D{h}$, and the 2-form  $\lie{X_h}\omega_R$. Imposing the requirement $\left.\lie{X_h}\omega_R \right|_{\D{h}}=0$ is a fairly straightforward matter now: it gives rise to polynomial expressions in the $(q,p)$-variables, the coefficients of which all have to vanish. Nevertheless, these computations are quite tedious so that assistance of Maple (or any other computer algebra package) is a great asset. We will not give a full account of the calculations involved, but merely indicate the order in which consecutive information is gathered, which will ultimately lead to the identification of admissible Hamiltonians of the form (\ref{exH}). The condition $\lie{X_h}\omega_R (X_h,\wt{R}(X_h))=0$ gives rise to a polynomial expression of degree 2 in the $p_i$; the coefficients are themselves polynomials in the $q^i$, with a degree which varies from 3 to 6. It turns out that the $p_1^2$-terms generate, among other conditions, a first-order differential equation for $c_1$ and algebraic relations which fix $a_3$ in terms of $a_1$ and $a_4$ in terms of $a_2$. Explicitly, we find that
\[
c_1(t) = c\,t^2 - t^{-1}, \quad a_3(t) = 3t^{-1}\sigma_1(t)a_1(t), \quad a_4(t) = \onethird t^{-1}\sigma_1(t)a_2(t),
\]
where $c$ is a constant. Further useful info comes from some of the monomials of degree 4 in the coefficients of $p_1$ and $p_2$. We can easily integrate equations for $a_1$, $a_2$, $\sigma_1$ and $c_2$, leading to
\[
a_1(t) = \alpha_1 t^5, \quad a_2(t) = \alpha_2 t^{9/2}, \quad \sigma_1 = s_1, \quad c_2(t) = - \onehalf t^{-1} + cs_1\,t^2,
\]
where $\alpha_1$, $\alpha_2$ and $s_1$ are arbitrary constants again. A little side remark here is that the above conclusions about $\sigma_1$ and $c_2$ follow in the first place from equations having an overall factor $\alpha_1$. The case $\alpha_1=0$ thus requires a separate analysis, but turns out to produce in the end a subcase of the solution we are going to present below. There are numerous other coefficients to be looked at, but they are all zero provided we finally impose that $s_1=1$. Having settled that, there are two more polynomial relations to be investigated, coming from the requirements $\lie{X_h}\omega_R (X_h,\wt{R}^2(X_h))=0$ and $\lie{X_h}\omega_R (\wt{R}(X_h),\wt{R}^2(X_h))=0$. It is an intriguing observation, however, that these are identically satisfied as a result of the conclusions we drew from the first condition. We thus arrive at the following class of separable Hamiltonians,
\begin{align}
H &= \onehalf p_1^2 + \onehalf t\,p_2^2 + 2t^3(t\,q_1^2 + q_2^2) + (c\,t^2 - t^{-1})p_1q_1 + ( c\,t^2 - \onehalf t^{-1})p_2q_2 \nonumber \\[1mm]
& \mbox{} \qquad + \alpha_1 t^5 q_1^3 + \alpha_2 t^{9/2} q_1^2q_2 + 3\alpha_1 t^4 q_1q_2^2 + \onethird \alpha_2 t^{7/2} q_2^3, \label{finalexH}
\end{align}
where $c$, $\alpha_1$ and $\alpha_2$ are arbitrary constants. The $R$-tensor which guarantees separability reads
\begin{align}
R &= t q_1 \left(\fpd{}{q_1}\otimes dq_1 + \fpd{}{q_2}\otimes dq_2\right) + q_2\left(\fpd{}{q_1}\otimes dq_2 + t\fpd{}{q_2}\otimes dq_1\right) \nonumber \\
& \mbox{} \qquad + (q_1^2 + \onehalf t^{-1}q_2^2) \fpd{}{q_1}\otimes dt +  \threehalf q_1q_2\fpd{}{q_2}\otimes dt. \label{finalexR}
\end{align}
Let us finally put our claim to an explicit test. As discussed before, the eigenvalues of $R$ suggest a change of coordinates, which reads
\[
Q_1 = t\,q_1 + \sqrt{t}\, q_2, \qquad Q_2 = t\,q_1 - \sqrt{t}\, q_2.
\]
The induced change of momenta can be written in the form
\[
p_1 = t(P_1 + P_2), \qquad p_2 = \sqrt{t}(P_1 - P_2).
\]
From (\ref{K}), we subsequently find that the Hamiltonian $K(t,Q,P)$ of the transformed system is given by
\begin{align}
K(t,Q,P) &= t^2\Big( P_1^2 + P_2^2 + Q_1^2 + Q_2^2 + c(P_1Q_1 + P_2Q_2) \nonumber \\
& \mbox{} \qquad + \onehalf \alpha_1(Q_1^3 + Q_2^3) + \onesix \alpha_2(Q_1^3 - Q_2^3)\Big). \label{finalexK}
\end{align}
It is clear that the Hamilton-Jacobi equation for $K$ can be solved by separation of variables indeed.

\subsubsection*{Acknowledgements} This work is part of the IRSES
project GEOMECH (nr. 246981) within the 7th European Community
Framework Programme. One of us (W.S.) further acknowledges support from the Czech Science Foundation under grant GACR 14-024765 ``Variations, Geometry and Physics".

\end{document}